# Geometric Interpretations and Applications of the Berger–Ebin and York $L^2$-Orthogonal Decompositions

**Sergey E. Stepanov**[1,2][*], **Irina I. Tsyganok**[2]

[1] Department of Mathematics, Russian Institute for Scientific
and Technical Information of the Russian Academy of Sciences,
20, Usievicha street, 125190 Moscow, Russia,
E-mail address: s.e.stepanov@mail.ru

[2] Department of Mathematics,
Finance University, 49-55, Leningradsky Prospect, 125468 Moscow, Russia,
E-mail address: i.i.tsyganok@mail.ru

**Abstract.** The Berger–Ebin and York $L^2$-orthogonal decompositions of the vector space of symmetric bilinear differential two-forms are fundamental tools in global Riemannian geometry. In this paper, we investigate the structure of Ricci tensors on compact Riemannian manifolds, with a particular focus on compact Ricci almost solitons, utilizing both the Berger–Ebin and York $L^2$-orthogonal decompositions. In addition, we explore applications of the York $L^2$-orthogonal decomposition to the theory of submanifolds and to the study of harmonic maps between Riemannian manifolds.

## 1. Introduction

The *Berger–Ebin* and *York* $L^2$-orthogonal decompositions of the vector space of symmetric bilinear differential two-forms (see [1], [2]) are foundational results in global Riemannian geometry and have been extensively studied in the literature (see, e.g., [3]). A key insight in Berger and Ebin's approach is the observation that, for any fixed Riemannian metric $g$ on a smooth compact manifold $M$, the space of differential symmetric two-form $C^\infty(S^2 M)$ admits a decomposition into two complementary closed subspaces that are orthogonal with respect to the $L^2$-inner product induced by $g$. This orthogonal splitting plays a central role in the study of Riemannian functionals (see [3, Ch. 4]), providing the infinitesimal counterpart of the slice theorem for the space of Riemannian metrics on $M$ (see [4]).

_________________________

[*] Corresponding author.

E-mail address: s.e.stepanov@mail.ru (S.E. Stepanov)

The *York decomposition* further refines this structure by expressing any differential symmetric two-form on a compact Riemannian manifold $(M, g)$ as a unique sum of a transverse-traceless component, a longitudinal (or pure gauge) component, and a pure trace component (see [1]). These components are determined via differential operators and are mutually orthogonal with respect to the $L^2$-inner scalar product. The covariant orthogonal decomposition of symmetric tensors is of particular interest in both mathematical and physical contexts, notably in general relativity (see, e.g., [5], [6], [8]) and in the classification theory of Riemannian metrics (see [1], [3], [7]).

The aim of this paper is to explore the *geometrization* of the Berger–Ebin and York decompositions through various significant applications in global Riemannian geometry. Specifically, we investigate the structure of Ricci tensors on compact Riemannian manifolds—focusing in particular on *Ricci almost solitons* — by employing these $L^2$-orthogonal decompositions. Furthermore, we study the decompositions of the *second fundamental form* of hypersurfaces in Riemannian manifolds of constant sectional curvature, illustrating their relevance to submanifold geometry. Finally, we consider applications of the York decomposition in the context of *harmonic maps* between Riemannian manifolds, demonstrating the utility of these decompositions in analytical and geometric settings.

## 2. Basic concepts and notations of Berger–Ebin and York $L^2$-orthogonal decompositions of symmetric bilinear differential 2-forms

In this preliminary section, we give a short account of some basic concepts and notations which will be used throughout the article.

We denote by $S^p M := S^p(T^*M)$ the vector bundle of differential symmetric bilinear $p$-forms or, in other words, of differential covariant symmetric $p$-tensors $(p \geq 1)$ on a compact Riemannian manifold $(M, g)$ and define the $L^2$-*inner scalar product* of symmetric bilinear differential $p$-forms $\varphi$ and $\phi$ on $(M, g)$ by the formula

$$\langle \varphi, \varphi' \rangle := \int_M g(\varphi, \phi) \, dvol_g$$

where $dvol_g$ being the volume element of $(M, g)$. Also $\delta^*: C^\infty(T^*M) \to C^\infty(S^2M)$ will be the first-order differential operator defined by the formula $\delta^*\theta := \frac{1}{2}L_\xi g$, where $L_\xi$ is the Lie derivative and $\xi = \theta^\#$ is the vector field dual (by $g$) to the 1-form (see [2]; [3, p. 117; 35; 514]). At the same time, we denote by the formula $\delta: C^\infty(S^2M) \to C^\infty(T^*M)$ the formal adjoint operator for $\delta^*$ which is called the *divergence of symmetric two-tensors* or, in other words, differential symmetric two-forms (see [2]; [3, p. 35]). In this case, we have $\langle \varphi, \delta^*\theta \rangle = \langle \delta\varphi, \theta \rangle$ for any $\varphi \in C^\infty(S^2M)$ and $\theta \in C^\infty(T^*M)$.

Since the symbol of the operator $\delta^*$ is injective (see [2]; [3, p. 117]), it follows from [2] and [3, Appendix] that the space $C^\infty(S^2M)$ admits the following direct sum decomposition

$$C^\infty(S^2M) = \operatorname{Im} \delta^* \oplus \delta^{-1}(0), \qquad (2.1)$$

which orthogonal for the $L^2$-inner scalar product (see [2]; [3, p. 118]). This is the *Berger–Ebin theorem* [2], which is a classical result in the Riemannian geometry in the large. By the Berger–Ebin theorem, an arbitrary symmetric bilinear form $\varphi \in C^\infty(S^2M)$ defined on a compact Riemannian manifold $(M,g)$ can be represented as the $L^2$-orthogonal sum $\varphi = \delta^*\theta + \varphi_0$, where $\varphi_0$ and $\theta$ are divergence free differential symmetric two-form and some differential one-form, respectively. Moreover, the equation holds $\langle \delta^*\theta, \varphi_0 \rangle = \langle \theta, \delta\varphi_0 \rangle = 0$ since $\delta\varphi_0 = 0$. Forthermore, the vector field $\xi = \theta^\#$ is defined up to a *Killing vector field* (see [3, p. 35]; [15, p. 43]). We recall here that a vector field $X$ on a Riemannian manifold $(M, g)$ is called an infinitesimal isometry transformation or Killing vector field if it generates a local one-parameter group of local isometric transformations (see [3, p. 40]). This means, that $L_X g = 0$ (see also [3, p. 40]; [15, p. 43]). Therefore, we have the following $\delta^*\theta = \frac{1}{2}L_\xi g = \frac{1}{2}L_{\xi+X} g$.

We will call formula (2.1) the *Berger–Ebin $L^2$-orthogonal decomposition for the vector space of differential symmetric* 2-*forms* on a compact Riemannian manifold.

In turn, *York's theorem* [1] is another well-known result and is also included in the monographs (see, e.g., [3, p. 130]). Namely, for any $n$-dimensional ($n \geq 3$) compact Riemannian manifold $(M, g)$ the decomposition holds (see [1])

$$C^{\infty}(S^2 M) = (\text{Im } \delta^* + C^{\infty}M \cdot g) \oplus \left(\delta^{-1}(0) \cap \text{trace}_g^{-1}(0)\right) \quad (2.2)$$

where both factors are infinite dimensional and orthogonal to each other with respect to the $L^2$ inner scalar product (see [3, p. 130]). Furthermore, the second factor $\delta^{-1}(0) \cap trace_g^{-1}(0)$ of (2.2) is the space of *TT-tensors*.

**Remark.** We recall that a symmetric divergence free and traceless covariant two-tensor is called *TT*-tensor (see, for example, [9]). The vector space of *TT*-tensors $\varphi^{TT}$ is defined by the condition (see, e.g., [10])

$$S^{TT}(M) := \{\varphi \in C^{\infty}(S^2 M) |\, \delta\,\varphi = 0, trace_g \varphi = 0\}.$$

As a consequence of a result of Bourguignon-Ebin-Marsden (see [3, p. 132]; [11]) the vector space $S^{TT}(M)$ is an infinite-dimensional vector space on any compact Riemannian manifold $(M, g)$. Such tensors are of fundamental importance in General Relativity (see, for example, [1]; [5]; [6] and [12]) and in Riemannian geometry (see, for instance, [3, p. 346-347]; [7] and [10]).

If we suppose $\varphi \in C^{\infty}(S^2 M)$, then York $L^2$-orthogonal decomposition formula (2.2) can be rewritten in the form

$$\varphi = \left(\frac{1}{2} L_\xi g + \lambda\, g\right) + \varphi^{TT} \quad (2.3)$$

for some $\xi \in C^{\infty}(TM)$, some $\varphi^{TT} \in S^{TT}(M)$ and some scalar function $\lambda \in C^{\infty}(M)$. Applying the operator $trace_g$ to both sides of (2.3), we obtain

$$trace_g \varphi = -\delta\theta + n\,\lambda, \quad (2.4)$$

where $\theta$ is the $g$-dual one-form of $\xi$. In this case, (2.3) can be rewritten in the form

$$\overset{\circ}{\varphi} = S\theta + \varphi^{TT} \quad (2.5)$$

where $\overset{\circ}{\varphi} = \varphi - \frac{1}{n}(trace_g \varphi)\, g$ is the traceless part of $\varphi$ and $S\theta := \delta^*\theta + \frac{1}{n}\delta\theta\, g$ denotes the *Cauchy-Ahlfors operator* $S: C^{\infty}(T^*M) \to C^{\infty}(S_0^2 M)$ actions on the vector space of one-form $C^{\infty}(T^*M)$ and with values in the vector space $C^{\infty}(S_0^2 M)$ of symmetric traceless bilinear differential forms (see, for example, [3]). It's obvious that $S$ annihilates the one-form $\theta$ such that $\theta^{\#} = \xi$ for a *conformal Killing vector field* $\xi$ on $(M, g)$, since the conformal Killing vector field $\xi$ obeys the equation $\delta^* \theta = -\frac{1}{n} \delta\theta \cdot g$ (see [15, p. 46]). Then the kernel of $S$ is the finite-dimensional vector space on a compact Riemannian manifold $(M, g)$. Therefore, for any $n$-dimensional $(n \geq 3)$ compact Riemannian manifold $(M, g)$ the decomposition holds

$$C^{\infty}(S_0^2 M) = \text{Im } S \oplus \left( \delta^{-1}(0) \cap \text{trace}_g^{-1}(0) \right), \qquad (2.6)$$

where both terms on the right-hand side of (2.2) are $L^2$-orthogonal to each other.

## 3. The Berger–Ebin $L^2$-orthogonal decomposition of the Ricci tensor and its applications

Let $(M, g)$ be an $n$-dimensional $(n \geq 3)$ compact Riemannian manifold and $Ric$ its Ricci tensor. Like the metric tensor $g$, the Ricci tensor $Ric$ assigns to each tangent space of the manifold a symmetric bilinear form (see [3, p. 43]). Therefore, from the Berger–Ebin $L^2$-orthogonal decomposition (2.1) of symmetric bilinear differential 2-forms we can deduce the $L^2$-orthogonal decomposition of the Ricci tensor as a symmetric bilinear form. Namely, we have the following

$$Ric = \delta^* \theta + \varphi_0 \qquad (3.1)$$

for some one-form $\theta \in C^{\infty}(T^*M)$ and some smooth divergence free symmetric two-form $\varphi_0 \in C^{\infty}(S^2 M)$. Applying the operator $\delta$ to both sides of (3.1), we obtain

$$\delta\delta^* \theta = -\frac{1}{2} ds \qquad (3.2)$$

since $\delta\varphi_0 = 0$ in accordance with the assumption made above. From equation (3.2) we deduce the integral formula

$$\int_M \left( L_\xi s \right) dv_g = - \langle \delta^*\theta, \delta^*\theta \rangle \geq 0 \qquad (3.3)$$

for $\xi = \theta^{\#}$. In this case, if $\int_M \left( L_\xi s \right) dv_g = 0$, then from (3,3) we deduce $\delta^*\theta := \frac{1}{2} L_\xi g = 0$ and hence $\xi$ is a Killing vector field. The opposite is also true. Furthermore, if $\xi$ is a Killing vector field, then from (3.2) we obtain that the scalar curvature $s$ of $(M, g)$ is constant. In addition, the opposite is also true. In conclusion, we recall that

We proved the following theorem.

**Theorem 1.** *Let $(M, g)$ be an $n$-dimensional $(n \geq 3)$ compact Riemannian manifold and (3.1) is the Berger–Ebin $L^2$-orthogonal decomposition for the Ricci tensor Ric of $(M, g)$ for some vector field $\xi \in C^\infty(TM)$ and some divergence free symmetric two-form $\varphi_0 \in C^\infty(S^2 M)$. Then the following propositions are equivalent:*

(i) *the scalar curvature $s$ of $(M, g)$ is constant;*

(ii) *the vector field $\xi$ is Killing;*

(iii) $\int_M (L_\xi s) \, dv_g = 0.$

Next, we recall the well-known Yamabe problem: Let $(M, g)$ be a compact Riemannian manifold, then there exists a positive and smooth function $f$ on $M$ such that the Riemannian metric $\bar{g} := f \cdot g$ has the constant scalar curvature $\bar{s}$. In turn, Richard Schoen provided in [9] an affirmative resolution to the problem in 1984. Thus, from the above it follows that if $(M, g)$ is a compact Riemannian manifold, then there exists a positive and smooth function $f$ on $M$ such that the Berger–Ebin $L^2$-orthogonal decomposition for the Ricci tensor $\overline{Ric}$ of the Riemannian metric $\bar{g} = f \cdot g$ can be rewritten in the form $\overline{Ric} = \bar{\varphi}_0$.

**Remark.** According to Theorem 1, we conclude that the presence of the metric $\bar{g} = f \cdot g$ (where $f$ is a positive and smooth function on $M$) with constant scalar curvature $\bar{s}$ on a compact Riemannian manifold $(M, g)$ is related to the existence (in general) of a non-zero Killing vector field on $(M, g)$.

## 4. The York $L^2$-orthogonal decomposition of the Ricci tensor and its applications

In turn, from the York $L^2$-orthogonal decomposition (2.2) we deduce the following $L^2$-orthogonal decomposition of the Ricci tensor:

$$Ric = (\delta^* \theta + \lambda g) + \varphi^{TT} \qquad (4.1)$$

for some one-form $\theta \in C^\infty(T^*M)$, some $TT$-tensor $\varphi^{TT} \in C^\infty(S^2 M)$ and some scalar function $\lambda \in C^\infty(M)$. Assume that $\xi = \theta^\#$, then it is obvious that $\delta^* \theta := \frac{1}{2} L_\xi g =$

$\frac{1}{2}L_{\xi+X}g$ for any Killing vector field $X$ and, therefore, the vector field $\xi$ is defined up to a *Killing vector field*. In addition, when $\xi$ is a Killing vector field the York's $L^2$-orthogonal decomposition for the Ricci tensor will be called *trivial*. In this case formula (3.4) can be rewritten as the form

$$Ric = \lambda g + \varphi^{TT}. \tag{4.2}$$

In this case $\lambda = \frac{s}{n}$ for the constant scalar curvature $s$ by the lemma of Schur (see [14, p. 14]; [3, p. 43]). Furthermore, if $\varphi^{TT} \equiv 0$, then the Ricci tensor $Ric$ satisfies the condition

$$Ric = \lambda g. \tag{4.3}$$

In this case $(M, g)$ is the *Einstein manifold* (see [3, p. 44]). Let us carry out the inverse reasoning. Namely, let $(M, g)$ be an $n$-dimensional ($n \geq 3$) compact Einstein manifold such that $Ric = \lambda g$. Then from (3.4) the equality follows $\delta^*\theta + \varphi^{TT} = 0$. Applying the operator $\delta$ to both sides of the last equality, we obtain $\delta\delta^*\theta = 0$. Note that $\ker \delta\delta^* = \ker \delta^*$ since $\langle \delta\delta^*\theta, \theta \rangle = \langle \delta^*\theta, \delta^*\theta \rangle$ for any $\theta \in C^\infty(T^*M)$. Therefore, the kernel of $\delta\delta^*$ (also as the kernel of $\delta^*$) is a vector space of Killing one-forms on a compact Riemannian manifold $(M, g)$. As a result from $\delta^*\theta + \varphi^{TT} = 0$ we deduce that $\delta^*\theta$ and $\varphi^{TT} = 0$.

**Proposition 1.** *An $n$-dimensional ($n \geq 3$) compact Riemannian manifold $(M, g)$ is Einstein if and only if, the vector field $\xi := \theta^\#$ is Killing and the TT-tensor $\varphi^{TT}$ is zero in the York decomposition $Ric = (\delta^*\theta + \lambda g) + \varphi^{TT}$ of its Ricci tensor.*

**Remark.** According to Proposition 1, we conclude that the definition of a compact Einstein manifold $(M, g)$ is related to the existence (in general) of a non-zero Killing vector field on $(M, g)$.

Next, we will need new concepts and related statements. Let $S$ be the *Cauchy-Ahlfors* operator actions on the one-form $\theta$ defined above (see [19]). The formal adjoint operator for $S$ is defined by the formula (see [17]; [18]; [19])

$$S^*\omega = 2\delta\omega$$

for an arbitrary traceless two-form $\omega \in C^\infty(S_0^2 M)$. Then the elliptic operator of the second kind $S^*S: C^\infty(T^*M) \to C^\infty(T^*M)$ is well known as the *Ahlfors Laplacian* (see [13]; [17]; [18]; [19] and [31]). Note that $\ker S^*S = \ker S$ since $\langle S^*S\theta, \theta \rangle = \langle S\theta, S\theta \rangle$ for any $\theta \in C^\infty(T^*M)$. Therefore, the kernel of $S^*S$ is also a finite dimensional vector space of conformal Killing one-forms on a compact Riemannian manifold $(M, g)$ (see [14]). Moreover, it is easy to verify the formula (see also [13]; [17]; [18] and [19])

$$S^*S\,\theta = 2\delta d\theta - 4Ric(\xi, \cdot) + \frac{4(n-1)}{n} d\delta\theta \qquad (4.4)$$

for a one-form $\theta$ and a vector field $\xi$ such that $\theta^\# = \xi$. Now, we are ready to formulate the lemma.

**Lemma 1.** *Let $(M, g)$ be an n-dimensional $(n \geq 3)$ compact Riemannian manifold and*

$$Ric = (\delta^*\theta + \lambda\, g) + \varphi^{TT}$$

*is the York decomposition for its Ricci tensor Ric, where $\xi$ is a vector field $\xi \in C^\infty(TM)$, $\varphi^{TT}$ is a TT-tensor $\varphi^{TT} \in C^\infty(S^2 M)$ and $\lambda$ is a scalar function $\lambda \in C^\infty(M)$. Then the following equation holds*

$$S^*S\,\theta = -\frac{n-2}{n} ds \qquad (4.5)$$

*for the one-form $\theta$ identified with $\xi$ by the Riemannian metric $g$, the scalar curvature $s$ of $(M, g)$ and the Ahlfors Laplacian $S^*S$.*

**Proof.** Applying the operator $trace_g$ to both sides of (4.4), we obtain

$$s = -\delta\theta + n\lambda \qquad (4.6)$$

where $\theta$ is the $g$-dual one-form of $\xi$. In this case, (4.4) can be rewritten in the form (see also (2.6))

$$\overset{\circ}{Ric} = S\theta + \varphi^{TT} \qquad (4.7)$$

where $\overset{\circ}{Ric} = Ric - \frac{1}{n} s\, g$ is the traceless part of the Ricci tensor $Ric$ and $S\theta = \frac{1}{2} L_\xi g + \frac{1}{n}\delta\theta\, g$ denotes the Cauchy-Ahlfors operator defined above. Next, applying $S^*$ to both sides of (3.7), we obtain (3.5), because $\delta\, Ric = -\frac{1}{2} ds$ by the Schur's lemma, which is a simple consequence of the "twice-contracted second Bianchi identity" (see [3, p. 43]). The proof of the lemma is complete.

Next, we will prove the Theorem 1.

**Theorem 2.** *Let $(M, g)$ be an n-dimensional $(n \geq 3)$ compact Riemannian manifold and (2.4) is the York $L^2$-decomposition for the Ricci tensor Ric of $(M, g)$, where $\xi$ is a vector field $\xi \in C^\infty(TM)$, $\varphi^{TT}$ is a TT-tensor $\varphi^{TT} \in C^\infty(S^2 M)$ and $\lambda$ is a scalar function $\lambda \in C^\infty(M)$. Then the following propositions are equivalent:*

(1) *the scalar curvature $s$ of $(M, g)$ is constant;*

(2) *the Ricci tensor Ric of $(M, g)$ has a trivial decomposition, i.e., $Ric = \frac{s}{n} g + \varphi^{TT}$;*

(3) $\int_M (L_\xi s)\, dv_g = 0.$

**Proof.** Let $(M, g)$ be an $n$-dimensional $(n \geq 3)$ compact Riemannian manifold and let (4.4) be the York decomposition for the Ricci tensor $Ric$ of $(M, g)$, then (4.7) holds. In turn, from (3.7) we deduce

$$-\frac{n-2}{n} \int_M (L_\xi s)\, dv_g = \langle S\theta, S\theta \rangle \quad (4.8)$$

for the one-form $\theta$ identified with $\xi$ by the Riemannian metric $g$ and the scalar curvature $s$ of $(M, g)$. Then we can conclude that the condition (1) and (3) are equivalent since (4.7) holds. Side by side, $\xi$ is a conformal Killing vector field. On the other hand, from (3) we deduce that $s$ is constant by the Schur's lemma (see [15, p. 14]; [3, p. 43]). In conclusion, we note that $\varphi_0 = \frac{s}{n} g + \varphi^{TT}$.

We can prove a corollary.

**Corollary 1.** *Let $(M, g)$ be an n-dimensional $(n \geq 3)$ compact Riemannian manifold and is the York $L^2$-decomposition for the Ricci tensor Ric of $(M, g)$, where $\xi$ is a vector field $\xi \in C^\infty(TM)$, $\varphi^{TT}$ is a TT-tensor $\varphi^{TT} \in C^\infty(S^2 M)$ and $\lambda$ is a scalar function $\lambda \in C^\infty(M)$. Then $\xi$ cannot satisfy the integral inequality $\int_M (L_\xi s)\, dv_g > 0$ for the scalar curvature $s$ of $(M, g)$.*

**Proof.** Let the conditions formulated above be fulfilled, then (4.8) holds. In this case the integral inequality $\int_M (L_\xi s) \, dv_g > 0$ contradicts the integral equality (4.8) since $\langle S\theta, S\theta \rangle \geq 0$.

**Remark.** For the case when $\xi \in C^\infty(TM)$ is a conformal Killing vector field an $n$-dimensional ($n \geq 3$) compact Riemannian manifold we obtain from (3.8) the equation $\int_M (L_\xi s) \, dv_g = 0$. Therefore, if $\int_M (L_\xi s) \, dv_g \neq 0$ for the vector field $\xi \in C^\infty(TM)$ from (2.4), then $\xi$ is not a conformal Killing vector field. We recall that the integral equality $\int_M (L_\xi s) \, dv_g = 0$ constitutes the main content of the well-known theorem in the classical article [20] on infinitesimal conformal transformations.

Taking into account the above, we will prove the theorem.

**Theorem 3.** *Let $(M, g)$ be a compact Riemannian manifold, then there exists a positive and smooth function $f$ on $M$ such that the York decomposition for the Ricci tensor $\overline{Ric}$ of the Riemannian metric $\bar{g} = f \cdot g$ is trivial, i.e.,*

$$\overline{Ric} = \frac{\bar{s}}{n} \bar{g} + \bar{\varphi}^{TT},$$

*where $\bar{\varphi}^{TT} \in C^\infty(S^2 M)$ is a TT-tensor on the Riemannian manifold $(M, \bar{g})$. Furthermore, the scalar curvature $\bar{s}$ of $(M, \bar{g})$ is constant.*

**Proof.** First, we note that (2.2) is a conformally invariant orthogonal decomposition (see [21]). In turn, Richard Schoen provided in [16] an affirmative resolution to the Yamabe problem in 1984: there exists a positive and smooth function $f$ on a compact manifold $(M, g)$ such that the Riemannian metric $\bar{g} := f \cdot g$ has the constant scalar curvature $\bar{s}$. Therefore, taking into account our Lemma and Theorem 2, we can conclude that Theorem 3 is true.

The Sampson Laplacian $\Delta_S : C^\infty(T^*M) \to C^\infty(T^*M)$ is defined by (see [22])

$$\Delta_S := 2\delta\delta^* - d\delta.$$

We recall that the vector field $\xi$ is an *infinitesimal harmonic transformation* in $(M, g)$ if the local one-parameter group infinitesimal transformations generated by the vector field $\xi$ is a group of harmonic transformations (see [22]). Furthermore, necessary and sufficient condition for a vector field $\xi$ on a Riemannian manifold to be infinitesimal harmonic transforma in $(M, g)$ is that $\theta \in Ker\, \Delta_S$, where $\theta$ is the $g$-dual one-form to $\xi$ (see [22]). We formulate a lemma needed to prove our next main results.

**Proposition 2**. *Let $(M, g)$ be an $n$-dimensional $(n \geq 3)$ compact Riemannian manifold and (2.4) is the York decomposition for the Ricci tensor Ric of $(M, g)$, where $\xi$ is a vector field $\xi \in C^\infty(TM)$, $\varphi^{TT}$ is a TT-tensor $\varphi^{TT} \in C^\infty(S^2M)$ and $\lambda$ is a scalar function $\lambda \in C^\infty(M)$. Then the equation holds*

$$\Delta_S \theta = -(n-2)d\lambda \tag{4.9}$$

*where $\theta$ is the $g$-dual one-form of $\xi$.*

**Proof.** We apply operator $trace_g$ to both sides of (14), we obtain equation (4.6), where $\theta$ is the $g$-dual one-form to the vector field $\xi$. Next, applying $d$ to both sides of (4.6), we obtain

$$ds = -d\delta\theta + n\, d\lambda. \tag{4.10}$$

In turn, applying the divergence $\delta$ to both sides of (4.1), we get

$$-ds = 2\delta\delta^*\theta - 2\, d\lambda. \tag{4.11}$$

Then using (4.10) from (4.11) we obtain $\Delta_S \theta = -(n-2)d\lambda$.

The following statement is obvious.

**Corollary 2**. *Let $(M, g)$ be an $n$-dimensional $(n \geq 3)$ compact Riemannian manifold and (2.4) is the York decomposition for the Ricci tensor Ric of $(M, g)$, where $\xi$ is a vector field $\xi \in C^\infty(TM)$, $\varphi^{TT}$ is a TT-tensor $\varphi^{TT} \in C^\infty(S^2M)$ and $\lambda$ is a scalar function $\lambda \in C^\infty(M)$. Then $\lambda = const$ if and only if $\xi$ is an infinitesimal harmonic transformation in $(M, g)$.*

## 5. Applications of the York $L^2$-orthogonal decomposition to the theory of Ricci almost solitons

One of the important components of the theory of Ricci flow are self-similar solutions called Ricci solitons (see [23, pp. 153-176]). Ricci solitons, which are a generalization of Einstein manifolds, have been studied more and more intensively in the last twenty years. This theory, besides being known after Perelman's proof of the Poincaré conjecture (for details see [24]), has a wide range of applications in differential geometry and theoretical physics. In turn, the study of Ricci almost solitons, which are a generalization of quasi-Einstein manifolds and Ricci solitons, was started by Pigola, Rigoli, Rimoldi, and Setty in [25].

An $n$-dimensional ($n \geq 2$) Riemannian manifold $(M, g)$ is called a *Ricci almost soliton*, if there exists a smooth complete vector field $\xi$ and exists a function $\lambda \in C^\infty(M)$ such that

$$Ric = \frac{1}{2} L_\xi g + \lambda\, g. \qquad (5.1)$$

Denote by $(M, g, \xi, \lambda)$ a Ricci almost soliton. For $\lambda = const$, it is a *Ricci soliton* (see [23]; [24] and etc.). Using Proposition 2, we can formulate a corollary about the compact Ricci almost soliton.

**Corollary 2**. *An $n$-dimensional ($n \geq 3$) compact Ricci almost soliton $(M, g, \xi, \lambda)$ is a Ricci soliton if and only if $\xi$ is an infinitesimal harmonic transformation in $(M, g)$.*

Note that the Killing vector field $\xi$ is an example of an infinitesimal harmonic transformation since it satisfies the following equations $\Delta_S \theta = 0$ and $\delta \theta = 0$ (see [22]). In particular, if $(M, g, \xi, \lambda)$ is an $n$-dimensional ($n \geq 3$) Ricci almost soliton with Killing vector field $\xi$, then $(M, g)$ is Einstein manifold, i.e., $Ric = \lambda\, g$ and $\lambda = \frac{1}{n} s$ for $s = const$ (see e.g., [15, p. 14] and [3, p. 43]). Moreover, when $\xi$ is a Killing vector field the almost Ricci soliton will be called *trivial*. Using the above, we can formulate a corollary of our Theorem 1.

**Corollary 3.** *Let $(M, g, \xi, \lambda)$ be an $n$-dimensional ($n \geq 3$) compact nontrivial Ricci almost soliton. Then the following propositions are equivalent:*

(1)  $\int_M (L_\xi s)\, dv_g = 0$ *for the scalar curvature $s$ of $(M, g)$;*

(2)  *the scalar curvature $s$ of $(M, g)$ is constant;*

(3)  *$(M, g)$ is an Einstein manifold, i.e., $Ric = \frac{s}{n} g$;*

(4)  *$(M, g)$ is isometric to a Euclidean sphere $\mathbb{S}^n$.*

**Proof.** In this corollary, the last statement needs to be proved. It in turn is a consequence of the theorem of Yano and Nagana (see [15, p. 57]): If a compact Einstein manifold

$(M, g)$ of dimension admits an infinitesimal nonisometric conformal transformation, then $(M, g)$ is isometric to a Euclidian sphere $\mathbb{S}^n$.

**Remark.** By Corollary 1, the vector field $\xi$ and the scalar curvature $s$ of $(M, g, \xi, \lambda)$ cannot satisfy the integral inequality $\int_M (L_\xi s)\, dv_g > 0$. Therefore, there is no compact Ricci almost soliton $(M, g, \xi, \lambda)$ satisfying the integral inequality above.

## 6. Applications of the York $L^2$-orthogonal decomposition to the submanifold theory of Riemannian manifolds

This section discusses Riemannian submanifolds. We will consider a smooth connected $n$-dimensional Riemannian manifold $(M, g)$ that is isometrically immersed into an $(n + 1)$-dimensional Riemannian manifold $(\bar{M}, \bar{g})$ for $n \geq 2$. Such a manifold is called a *hypersurface*. We denote by $\nabla$ and $\bar{\nabla}$ the Levi-Civita connection on $(M, g)$ and $(\bar{M}, \bar{g})$, respectively. In this case, the well-known Gauss formula says that

$$\nabla_X Y = \bar{\nabla}_X Y - g(A_g X, Y) N,$$

where $X, Y$ are any vector fields tangent to $(M, g)$, the vector field $N$ is the global unite vector field normal to $(M, g)$ and $A_g$ stands the *shape operator* of $(M, g)$. The *mean curvature* of $(M, g)$ is given by $H = trace A_g$ (see [3, p. 38]). Then $(M, g)$ is minimal if and only if $H \equiv 0$ (see [31, p. 39]). The second fundamental form of a hypersurface $(M, g)$ is defined by the identity $\varphi(X, Y) := g(A_g X, Y)$, where $X, Y$ are any vector fields tangent to $(M, g)$. If $(\bar{M}, \bar{g})$ is a Riemannian manifold of constant sectional curvature, then the second fundamental form $\varphi$ of $(M, g)$ satisfies the *Codazzi equations* (see [3, p. 436] and [26, p. 350]):

$$(\nabla_X \varphi)(Y, Z) = (\nabla_Y \varphi)(X, Z), \tag{6.1}$$

where $X, Y$ and $Z$ are any vector fields tangent to $(M, g)$. At the same time, any $\varphi \in C^\infty(S^2 M)$ satisfying Codazzi equations (1) is called the *Codazzi tensor* (see [3, p. 435]). If $(\bar{M}, \bar{g})$ is a Riemannian manifold of constant sectional curvature, then from Codazzi equations (1) we deduce

$$\delta \varphi = -\frac{1}{n}\, d\, H. \tag{6.2}$$

From the Berger–Ebin $L^2$-orthogonal decomposition (2.1) we deduce the decomposition for the second fundamental form $\varphi$ of $(M, g)$

$$\varphi = \delta^*\theta + \varphi_0 \tag{6.3}$$

for some one-form $\theta \in C^\infty(T^*M)$ and some smooth divergence free symmetric two-form $\varphi_0 \in C^\infty(S^2 M)$. Applying the operator $\delta$ to both sides of (6.3), we obtain

$$\delta\delta^*\theta = -\frac{1}{n} dH \tag{6.4}$$

since $\delta\varphi_0 = 0$ in accordance with the assumption made above. From equation (6.4) we deduce the integral formula

$$\int_M (L_\xi H) \, dv_g = -\langle \delta^*\theta, \delta^*\theta \rangle \geq 0 \tag{6.5}$$

for $\xi = \theta^\#$. In this case, if $\int_M (L_\xi H) \, dv_g = 0$, then from (6,5) we deduce $\delta^*\theta := \frac{1}{2} L_\xi g = 0$ and hence $\xi$ is a Killing vector field (see [15, p. 43]). The opposite is also true. Furthermore, if $\xi$ is a Killing vector field, then from (6.4) we obtain that mean curvature of $(M, g)$ is constant. In addition, the opposite is also true. We proved the following theorem.

**Proposition 3.** *Let $(M, g)$ be an $n$-dimensional ($n \geq 3$) compact Riemannian manifold that is isometrically immersed into an $(n+1)$-dimensional Riemannian manifold $(\bar{M}, \bar{g})$ with constant section curvature and (3.1) be the Berger–Ebin $L^2$-orthogonal decomposition for the second fundamental form $\varphi$ of $(M, g)$ for some vector field $\xi \in C^\infty(TM)$ and some divergence free symmetric two-form $\varphi_0 \in C^\infty(S^2 M)$. Then the following propositions are equivalent:*

(i) *the mean curvature $H$ of $(M, g)$ is constant;*

(ii) *the vector field $\xi$ is Killing;*

(iii) $\int_M (L_\xi H) \, dv_g = 0.$

**Remark.** According to Proposition 2, we conclude that the definition of a compact hypersurface $(M, g) \subset (\bar{M}, \bar{g})$ with constant non-zero mean curvature is related to the existence (in general) of a non-zero Killing vector field on $(M, g)$.

Let $(M, g)$ be a compact hypersurface with constant mean curvature in a Riemannian manifold $(\bar{M}, \bar{g})$ with constant sectional curvature and $\dim \bar{M} \geq 4$, then it follows that the second fundamental form $\varphi$ of $(M, g)$ is a *Codazzi tensor* with constant trace and zero divergence. At the same time, it is well-known that every Codazzi tensor $\varphi$ with constant trace on a compact Riemannian manifold $(M, g)$ with non-negative sectional curvature $sec$ is parallel. If, moreover, $sec > 0$ at some point, then $\varphi$ is a constant multiple of $g$ (see [3, p. 436]). At the same time, we recall that a Riemannian submanifold $(M, g)$ is said to be *parallel* if its second fundamental form $\varphi$ is parallel, that is $\nabla \varphi \equiv 0$ (see [26, Chapter 8]). On the other hand, a submanifold $(M, g)$ of a Riemannian manifold $(\bar{M}, \bar{g})$ is called *totally umbilical* if its second fundamental form $\varphi$ is proportional to its first metric tensor $g$ (see [3, p. 39]; [26, Chapter 12]). Using the definition above, we can formulate a statement.

**Corollary 4.** *Let $(M, g)$ be an n-dimensional $(n \geq 3)$ compact hypersurface with constant mean curvature in a Riemannian manifold $(\bar{M}, \bar{g})$ with constant sectional curvature. If the sectional curvature $sec$ of $(M, g)$ is non-negative, then it is a parallel submanifold. If, moreover, $sec > 0$ at some point, then $(M, g)$ is totally umbilical.*

**Remark.** In particular, if $(M, g)$ is a compact minimal hypersurface in a Euclidian sphere $\mathbb{S}^{n+1}$, and $(M, g)$ is of strictly positive sectional curvature, then it has to be an equator of $\mathbb{S}^{n+1}$ (see [2]).

On the other hand, the second fundamental form $\varphi \in C^\infty(S^2 M)$ has the York decomposition

$$\varphi = \left(\frac{1}{2} L_\xi g + \lambda\, g\right) + \varphi^{TT} \qquad (6.6)$$

for some vector field $\xi \in C^\infty TM$, $TT$-tensor $\varphi^{TT} \in C^\infty(S^2 M)$ and scalar function $\lambda \in C^\infty M$. If we applying the operator $trace_g$ to both sides of (6.6), we obtain

$$n\, H = -\delta\theta + n\, \lambda \qquad (6.7)$$

where $\theta^\# = \xi$ (see [31, p. 35]). In this case, (6.7) can be rewritten in the form

$$\overset{\circ}{\varphi} = S\theta + \varphi^{TT} \qquad (6.8)$$

where $\mathring{\varphi} = \varphi - H\,g$ is the traceless part of the second fundamental form $\varphi$ and $S\theta = \frac{1}{2}L_\xi g + \frac{1}{n}\delta\theta\, g$ is the Cauchy-Ahlfors operator. Next, applying $S^* := \delta$ to both sides of (6.7), we obtain

$$S^*S\,\theta = \delta\mathring{\varphi} \qquad (6.9)$$

for the Ahlfors Laplacian $S^*S$ (see [19]). As a result, we obtain

$$S^*S\,\theta = -\frac{1}{n}dH. \qquad (6.10)$$

Therefore, if the mean curvature $H$ of $(M, g)$ is constant, then $S^*S\,\theta = 0$. Therefore, $S\,\theta = 0$ since $\langle S^*S\,\theta, \theta\rangle = \langle S\,\theta, S\theta\rangle \geq 0$ (see [19]). In this case, from (6.8) we obtain $\varphi = H\,g + \varphi^{TT}$. The converse is also true.

**Theorem 4.** *Let $(M, g)$ be an $n$-dimensional $(n \geq 3)$ compact hypersurface in a Riemannian manifold $(\bar{M}, \bar{g})$ with constant sectional curvature. Then the mean curvature $H$ of $(M, g)$ is constant if and only if the second fundamental form $\varphi$ of $(M, g)$ admits the following $L^2$-orthogonal decomposition*

$$\varphi = H\,g + \varphi^{TT} \qquad (6.11)$$

*for some TT-tensor $\varphi^{TT}$. If, moreover, $(M, g)$ is a compact minimal hypersurface, then $\varphi$ is a TT-tensor $\varphi^{TT}$.*

From equation (6.11) and Codazzi equations, it follows that $\varphi^{TT}$ is a traceless Codazzi tensor, meaning that $\varphi^{TT}$ is a *harmonic form* (see [21, p. 350]). At the same time, we proved in [10] that any such $TT$-tensor must be the zero tensor on a compact Riemannian manifold $(M, g)$ with quasi-positive sectional curvature, i.e., it is semi-negative-definite everywhere and negative-definite at some point. In particular, if $(M, g)$ is a hypersurface in a Riemannian manifold $(\bar{M}, \bar{g})$ with constant sectional curvature, then it is totally umbilical with constant mean curvature. In this case $(M, g)$ is a Riemannian manifold with constant sectional curvature. Moreover, this curvature must be positive, given that we assumed that $(M, g)$ to be a manifold of quasi-positive sectional curvature. In our case, $(M, g)$ must be a spherical space form since it is a compact manifold of positive sectional curvature (see [26, p. 222]). In particular, if $(M, g)$ a simply connected manifold, then it is a Euclidean sphere $\mathbb{S}^n$. Using this fact, along with the results from Theorem 3, we arrive at the following conclusion:

**Corollary 5**. *Let $(\bar{M}, \bar{g})$ be an $n$-dimensional Riemannian manifold with constant sectional curvature, where $n \geq 4$, and let $(M, g)$ be an isometrically immersed hypersurfaces $(M, g) \to (\bar{M}, \bar{g})$ with constant mean curvature. If $(M, g)$ has quasi-positive sectional curvature, then $(M, g)$ is a spherical space form. In addition, if $(M, g)$ is a simply connected manifold, then it is a Euclidean sphere $\mathbb{S}^n$.*

## 7. Applications of the Berger–Ebin $L^2$-orthogonal decomposition to the theory of harmonic maps

Let $(M, g)$ be a compact oriented Riemannian manifold of dimension $n \geq 2$. Along with $(M, g)$, consider an $m$-dimensional ($m \geq 2$) Riemannian manifold $(\bar{M}, \bar{g})$ and a smooth $C^\infty$-mapping $f : (M, g) \to (\bar{M}, \bar{g})$. Let $f_*$ denote the differential of $f$. The differential $f_* : TM \to T\bar{M}$ is a $C^\infty$-section of the tensor bundle $T^*M \otimes f^*T\bar{M}$, where $f^*$ is the transpose map of $f_*$, and $f^*T\bar{M}$ the bundle with base $M$, fiber $T_{f(x)}\bar{M}$ over $x \in M$, and with the metric $\bar{g}'$ and connection $\bar{\nabla}'$ induced from $(\bar{M}, \bar{g})$ (see [27]).

The Riemannian structures $g$ and $\bar{g}$ define a metric $G$ on the fibers of the bundle $T^*M \otimes f^*T\bar{M}$. The connections $\nabla$ and $\bar{\nabla}'$ induce a connection $D$ in the bundle $T^*M \otimes f^*T\bar{M}$ such that (see [27, Theorem 1.4.1]; [28, formula (2.2)])

$$(Df_*)(X, Y) = \bar{\nabla}'_X f_* Y - f_* \nabla_X Y \qquad (7.1)$$

for any $X, Y \in C^\infty TM$.

It is known (see [29]) that a *harmonic map* $f : (M, g) \to (\bar{M}, \bar{g})$ is defined as an extremum of the Energy functional $E(f) = \frac{1}{2}\int_M (trace_g (f^*\bar{g})) \, dv_g$, and $f$ is a harmonic map if and only if it satisfies the following *Euler–Lagrange equation*: $trace_g(Df_*) = 0$, where $\tau(f) = trace_g(Df_*)$ is called the *tension field*. The following lemma is true (see [30]).

**Lemma 2**. *Let $f : (M, g) \to (\bar{M}, \bar{g})$ be a smooth submersion or diffeomorphism of a compact Riemannian manifold $(M, g)$ of dimension $n \geq 3$ onto an $m$-dimensional*

$(n \geq m)$ *Riemannian manifold* $(\bar{M}, \bar{g})$. *Then* $f$ *is a harmonic map if and only if the tensor field* $g^* = f^*\bar{g}$ *satisfy the equation* $\delta\left(g^* - \frac{1}{2}(trace_g\, g^*)\, g\right) = 0$.

**Proof.** The covariant derivative $\nabla g^*$, in accordance with (7.1), has the form (see also [30, formulas (2.2) and (2.3)])

$$(\nabla_Z g^*)(X,Y) = \bar{g}'\big((Df_*)(Z,X), f_*Y\big) + \bar{g}'\big((Df_*)(Z,Y), f_*X\big)$$

for any $X, Y, Z \in C^\infty TM$. It follows directly from this formula that

$$(\delta g^*)(Y) = -\bar{g}'(\tau(f), f_*Y) - \bar{g}'\big((Df_*)\, X_i, f_*X_i\big);$$
$$(\nabla\, trace_g\, g^*)(Y) = 2\, \bar{g}'\big((Df_*)\, X_i, f_*X_i\big)$$

for an arbitrary local orthonormal basis $X_1, \ldots, X_n$ of vector fields on $(M, g)$. From the above equalities, we deduce the following equality:

$$\delta\left(g^* - \frac{1}{2}(trace_g\, g^*)\, g\right)(Y) = -\bar{g}'(\tau(f), f_*Y). \tag{7.2}$$

Since $f$ is a harmonic submersion, by virtue of (7.2), we have $\delta\left(g^* - \frac{1}{2}(trace_g\, g^*)\, g\right) = 0$. To prove the converse statement, it is necessary to show that the equality $\bar{g}'(\tau(f), f_*Y) = 0$ implies the vanishing of the *tension field* $\tau(f) = trace_g(Df_*)$. This takes place when $Rank\, [(f_*)_x] \geq m$ for an arbitrary point $x \in M$, in particular, for a submersion. In conclusion of the proof, we note that we have already considered the case of a harmonic diffeomorphism $f: (M, g) \to (\bar{M}, \bar{g})$ in [31] and proved a similar proposition.

From the Berger–Ebin $L^2$-orthogonal decomposition (2.1) of symmetric bilinear differential two-forms we can deduce the $L^2$-orthogonal decomposition of the tensor $g^* = f^*\bar{g}$ as a symmetric bilinear form. Namely, we have the following decomposition:

$$g^* = \delta^*\theta + \varphi_0 \tag{7.3}$$

for some one-form $\theta \in C^\infty(T^*M)$ and some smooth divergence free symmetric two-form $\varphi_0 \in C^\infty(S^2 M)$. Applying the operator $\delta$ to both sides of (7.3), we obtain $\delta\delta^*\theta = \delta g^*$ and hence $2\, \delta g^* = \Delta_S \theta + d\delta\, \theta$. The last equation can be rewritten as

$$2\, \delta\left(g^* - \frac{1}{2}(trace_g g^*)\, g\right) = \Delta_S \theta + d\left(\delta\, \theta + (trace_g g^*)\right), \tag{7.4}$$

where, according to (7.2), the left-hand side of equality (7.4) is equal to $-2\,\bar{g}'(\tau(f), f_*Y)$. Therefore, we have the equation

$$-2\,\bar{g}'(\tau(f), f_*X) = (\Delta_S\theta)(X) + d\left(\delta\,\theta + (trace_g g^*)\right)(X) \quad (7.5)$$

for any $X \in C^\infty TM$. Note that $\Delta_S\theta = 0$ characterizes the vector field $\xi = \theta^\#$ as an infinitesimal harmonic transformation. Now the following theorem follows from the Lemma 2 proved above.

**Theorem 5.** *Let $f : (M, g) \to (\bar{M}, \bar{g})$ be a smooth submersion or diffeomorphism of a compact Riemannian manifold $(M, g)$ of dimension $n \geq 3$ onto an m-dimensional $(n \geq m)$ Riemannian manifold $(\bar{M}, \bar{g})$. Let $g^* = \delta^*\theta + \varphi_0$ be the Berger–Ebin $L^2$-orthogonal decomposition of differential symmetric two-form $g^* = f^*\bar{g}$ for some one-form $\theta \in C^\infty(T^*M)$ and some smooth divergence free symmetric two-form $\varphi_0 \in C^\infty(S^2M)$. If any two of the following three statements hold, then the third statement also holds:*

(i) *$f$ is a harmonic map;*

(ii) *$\xi = \theta^\#$ is an infinitesimal harmonic transformation;*

(iii) *$div\,\xi = trace_g g^* + C$ for an arbitrary constant $C$.*

As an immediate consequence of the theorem, we have the following corollary.

**Corollary 7.** *Let $f : (M, g) \to (\bar{M}, \bar{g})$ be a harmonic map between a compact Riemannian manifold $(M, g)$ of dimension $n \geq 2$ and an m-dimensional $(n \geq 2)$ Riemannian manifold $(\bar{M}, \bar{g})$. Let $g^* = \delta^*\theta + \varphi_0$ be the Berger–Ebin $L^2$-orthogonal decomposition of $g^* = f^*\bar{g}$, where $\theta \in C^\infty(T^*M)$ and $\varphi_0 \in C^\infty(S^2M)$ is a divergence free two-form. If $\xi = \theta^\#$ is an infinitesimal harmonic transformation, then the Energy functional of $f$ has the form $E(f) = C\,Vol(M, g)$ for an arbitrary constant $C \geq 0$.*

**Conflict of interest.** The authors declare that they have no conflict of interest.

**Funding Declaration.** The authors declare that no funds, grants, or other support were received during the preparation of this manuscript.